\documentclass[12pt]{amsart}
\usepackage{amsmath,amssymb,enumerate} 
\newtheorem{thm}{Theorem}[section]
\newtheorem{prop}[thm]{Proposition}
\newtheorem{lem}[thm]{Lemma}

\theoremstyle{definition}

\theoremstyle{remark}

\numberwithin{equation}{section}
\begin{document}
\title[Pulse solutions of a coupled NLS system]{
Solitary pulse solutions of a coupled nonlinear Schr\"{o}dinger system arising in optics}

\author[Sharad Silwal]{Sharad Silwal}
\address{Jefferson College of Health Sciences,
101 Elm Ave SE, Roanoke, VA 24013, USA} \email{sdsilwal@jchs.edu}

\begin{abstract}
We prove the existence of travelling-wave solutions for a system of coupled
nonlinear Schr\"{o}dinger equations arising in nonlinear optics. Such a system describes second-harmonic
generation in optical materials with $\chi^{(2)}$ nonlinearity.
To prove the existence of travelling waves, we employ the method of concentration compactness to
prove the relative compactness of minimizing sequences of the associated variational problem.
\end{abstract}

 \maketitle

\section{Introduction}

The coupled nonlinear Schr\"{o}dinger system we are studying in this paper has applications in nonlinear optics. Telecommunications and computer technology are prime examples where modern nonlinear optics plays a major role. The system we are considering has $\chi^{(2)}$ nonlinearity, which arises in an optical process where a phenomenon of frequency doubling is exhibited due to interactions beteween certain nonlinear materials. We refer the reader to \cite{[Men]} for a detailed description of $\chi^{(2)}$ nonlinearity, which, in fact, has gained the focus of many mathematicians and physicists in recent years, in an attempt to increase the speed and efficiency of optical fibres for data transmission. The reader may consult \cite{[Fer], [Men], [Yew]} and the references therein for more information about the physics and the engineering applications of the system studied in this paper.

In this paper, we consider a system of two-coupled nonlinear
Schr\"{o}di-nger equations in the form
\begin{equation}\label{nlseq}
\left\{
\begin{aligned}
& i\frac{\partial W}{\partial t} +\frac{\partial^2W}{\partial x^2}-\alpha W+\overline{W}V=0, \\
& i\mu \frac{\partial V}{\partial t} +\frac{\partial^2V}{\partial x^2} -\beta V+\frac{1}{2}W^2=0,
\end{aligned}
\right.
\end{equation}
where $W$ and $V$ are complex-valued functions, $\alpha$ and $\beta$ are real numbers, and $\mu>0.$
System \eqref{nlseq} is obtained from the basic $\chi^{(2)}$ second-harmonic generation equations (SHG) of type I (see \cite{[Men]}).
Physically, the complex functions $W$ and $V$ represent packets of
amplitudes of the first and second harmonics of an optical wave respectively.
The constant $\mu$ measures the ratios of the dispersions/diffractions.

\smallskip

In this paper, we are interested in the existence problem of travelling solitary wave solutions of \eqref{nlseq}. The travelling waves
we are interested in are of the form
\begin{equation}\label{SOdef}
(W(x,t), V(x,t))=(e^{i\sigma t}\Phi(x), e^{2i\sigma t}\Psi(x)),
\end{equation}
where $\Phi,\Psi:\mathbb{R}\to \mathbb{R}.$ We specify the boundary conditions $\Phi,\Psi\to 0$ as $x\to \pm\infty$ and call these solutions as pulses.
We show that there exists a nontrivial smooth solution with exponential decay at infinity for
the system
\begin{equation}\label{auxsys}
\left\{
\begin{aligned}
& \Phi^{\prime \prime}-\alpha_0 \Phi + \Phi \Psi=0,\\
& \Psi^{\prime \prime}-\beta_0 \Psi +\frac{1}{2}\Phi^2=0,
\end{aligned}
\right.
\end{equation}
with $\alpha_0=\alpha+\sigma>0$ and $\beta_0=\beta+2\mu \sigma>0.$
It is well-known that when $\alpha_0=\pm 1 = \beta_0$ and $\Phi=\pm \sqrt{2}\Psi,$
system \eqref{auxsys} has the explicit solutions of the form
\begin{equation*}
\Phi(x)=\pm\frac{3}{\sqrt{2}}\textrm{sech}^2(x/2),\ \ \ \Psi(x)=\pm \frac{3}{2}\textrm{sech}^2(x/2).
\end{equation*}
From the mathematical point of view, the existence of at
least one solution homoclinic to the origin (corresponding physically to a
pulse) for \eqref{auxsys} was
proved in \cite{[Yew]} for all $\beta_0>0$ and $\alpha_0=1.$ Their method is of variational nature and uses the mountain pass theorem,
along with some convergence arguments. Localized solutions were characterized as critical points which are local minima of an energy functional.
In \cite{[ABha]}, using the concentration compactness technique, the existence of pulses was established for a more general version of system \eqref{auxsys}.
Their method uses variational characterizations of pulses subject to two independent constraints. Solutions were characterized as critical points which are global minima of an energy functional, with $\alpha_0, \beta_0$ appearing as Lagranges multipliers
associated with the constraints.
To establish the compactness result for every minimizing sequence,
they showed the subadditivity condition through a new idea based on rearrangement inequalities.
In \cite{[Yew2]}, a description of the profile of solutions for \eqref{auxsys} was given by using
the framework of homoclinic bifurcation theory. In \cite{[AngF]}, the existence of periodic pulses
of \eqref{auxsys} as well as the stability and instability
of such solutions were studied. See also \cite{[SanDCDS]} where the author obtained 
positive solutions $(\Phi,\Psi)$ of \eqref{auxsys} satisfying the additional condition $\|\Phi\|_{2}^2=a$ and $\|\Psi\|_{2}^2=b$ for any $a>0$ and $b>0.$

\smallskip

In the present paper, we study a different variational problem than those
considered in \cite{[ABha],[Yew]} and establish the existence of travelling pulse solutions to \eqref{auxsys} for all $\alpha_0>0$ and $\beta_0>0.$
Our method exhibits a new one-parameter
family of travelling solitary-wave solutions than those obtained in \cite{[ABha],[Yew]}.
The key tool in our analysis is the concentration
compactness lemma of \cite{[L1]}. To rule out the dichotomy case while applying the
concentration compactness lemma, we use an argument developed in \cite{[LEV]}.
Similar techniques have been used previously in \cite{[ChenNg]} to prove the existence of travelling-wave
solutions to Boussinesq systems
and in \cite{[SB1]} to study solitary waves for an equation of short and long dispersive waves arising in two-layer fluids.

\smallskip

\noindent We now state our main results. Our existence result is proved by using a variational approach.
Precisely, let the functional $I:H^1(\mathbb{R})\times H^1(\mathbb{R})\to \mathbb{R}$ be defined by
\begin{equation}\label{Idef}
I(f,g)=\int_{\mathbb{R}}\left((f^\prime)^2+(g^\prime)^2+\alpha_0 f^2+\beta_0 g^2\right)\ dx,
\end{equation}
where $H^1(\mathbb{R})$ denotes the $L^2-$based Sobolev space of first order.
For any $\lambda>0$, we consider the minimization problem
\begin{equation}\label{varpro}
M_\lambda = \inf\left\{I(f,g):(f,g)\in H^1(\mathbb{R})\times H^1(\mathbb{R}),\ \int_{\mathbb{R}}f^2g\ dx=\lambda \right\}.
\end{equation}

\smallskip

\noindent The following is our existence result.

\begin{thm}
\label{existence}
For any $\lambda>0,$ define
\begin{equation}
\mathcal{S}_\lambda=\left\{(f,g) \in H^1(\mathbb{R})\times H^1(\mathbb{R}): I(f,g)=M_\lambda,\ \int_{\mathbb{R}}f^2 g \ dx=\lambda \right\}.
\end{equation}
Then the following statements hold.

\smallskip

\noindent (i)
The infimum $M_\lambda$ defined in \eqref{varpro} satisfies $0<M_\lambda<\infty.$

\smallskip

\noindent (ii) For every sequence $\{(f_{n},g_{n})\}$ in $H^1(\mathbb{R})\times H^1(\mathbb{R})$ such that
\begin{equation}\label{minseq}
\lim_{n\to \infty}\int_{\mathbb{R}}f_{n}^2g_n=\lambda\ \textrm{ and }\ \lim_{n\to
\infty }I(f_n,g_n)=M_\lambda,
\end{equation}
there exists a sequence of real numbers $\{y_{k}\}$ such that $\{(f_{n_k}(\cdot
+y_{k}),g_{n_k}(\cdot +y_{k})\}$ converges strongly in $H^1(\mathbb{R})\times H^1(\mathbb{R})$ to
some $(\phi,\psi)$ in $\mathcal{S}_\lambda$.  In particular, the set
$\mathcal{S}_\lambda$ is non-empty.

\smallskip

\noindent (iii) Each function $(\phi,\psi)$ in $\mathcal{S}_\lambda,$ after multiplying by a constant, is a solution of
\eqref{auxsys}, and hence when substituted into \eqref{SOdef} gives one-parameter
family of travelling-wave solutions to the system \eqref{nlseq}.

\smallskip

\noindent (iv) If $(\phi,\psi)\in\mathcal{S}_\lambda,$ then $\phi, \psi\in H^\infty(\mathbb{R})$.
and $\phi, \psi$ decay exponentially at infinity.
 \end{thm}

We now provide some notations that will be used throughout the paper.

\smallskip

\noindent \textit{Notation}.
For $1\leq p\leq \infty ,$ we denote by $L^{p}=L^{p}(\mathbb{R})$ the space
of all measurable functions $f$ on $\mathbb{R}$ for which the norm $%
\left\vert f\right\vert _{p}$ is finite, where%
\begin{equation*}
\left\vert f\right\vert _{p}=\left( \int_{-\infty }^{\infty }\left\vert
f\right\vert ^{p}dx\right) ^{1/p}\textrm{ \ for }1\leq p<\infty
\end{equation*}
and $\left\vert f\right\vert _{\infty }$ is the essential supremum of $%
\left\vert f\right\vert $ on $\mathbb{R}.$ For $s\geq 0,$ we denote by $H_{%
\mathbb{C}}^{s}(\mathbb{R})$ the Sobolev space of all complex-valued
functions $f$ in $L^{2}$ for which the norm%
\begin{equation*}
\left\Vert f\right\Vert _{s}=\left( \int_{-\infty }^{\infty }\left(
1+\left\vert \xi \right\vert ^{2}\right) ^{s}\left\vert \widehat{f}(\xi
)\right\vert ^{2}d\xi \right) ^{1/2}<\infty .
\end{equation*}
We will always view $H_{\mathbb{C}}^{s}(\mathbb{R})$ as a vector space over
the reals, with inner product given by%
\begin{equation*}
\left\langle f_{1},f_{2}\right\rangle =\textrm{Re}\int_{-\infty }^{\infty
}\left( 1+\left\vert \xi \right\vert ^{2}\right) ^{s}\widehat{f_{1}}%
\overline{\widehat{f_{2}}}\ d\xi .
\end{equation*}
The space of all real-valued functions $f$ in $H_{\mathbb{C}}^{s}(\mathbb{R})
$ will be denoted by $H^{s}(\mathbb{R}).$ In particular, we use $\left\Vert
f\right\Vert $ to denote the $L^{2}$ or $H^{0}(\mathbb{R})$ norm of a
function $f.$ We
define the space $X$ to be $H^{1}(\mathbb{R})\times H^{1}(%
\mathbb{R})$ provided with the product norm.
The letter $C$ will frequently be used to denote various constants whose
actual value is not important for our purposes.

 \section{Proof of the existence result}

\noindent In this section we provide the proof of our existence result.
Our strategy here is to employ the concentration
compactness method to solve the variational problem \eqref{varpro}. In this approach one
obtains travelling solitary waves as global minimizers of the problem. When the method works, it shows not only
that global minimizers exist, but also that every minimizing sequence is relatively
compact up to translation. In the last couple of decades, this method was applied by
many authors to prove the existence of solutions for a variety of
dispersive evolution equations (see, for example \cite{[A], [ABha], [SanDCDS], [ChenNg], [LEV], [Shar]}).

\smallskip

\noindent As usual in this method, we take
any minimizing sequence $\{(f_{n},g_{n})\}$ for $M_\lambda$ and consider a sequence of nondecreasing functions
$P_{n}:[0,\infty)\to [0,\lambda]$ defined by
\begin{equation*}
P_n(\omega)=\sup_{y\in \mathbb{R}}\int_{y-\omega}^{y+\omega}\rho _{n}(x)\ dx.
\end{equation*}
where $\rho _{n}(x)=(f_{n}'(x))^2 + (g_{n}'(x))^2+f_n^2(x)+g_n^2(x).$ As $P_{n}(\omega)$ is a uniformly bounded sequence of
nondecreasing functions in $\omega,$ one can show that it has a subsequence, which we will again denote by $P_n$,
that converges pointwise to a nondecreasing limit function $P(\omega):[0,\infty)\to [0,\lambda]$. Define
\begin{equation}
\label{defgamma} \gamma =\lim_{\omega\to \infty }P(\omega).
\end{equation}
Then $\gamma$ satisfies $0\leq \gamma \leq \lambda$. From Concentration-Compactness Lemma of P.~L.~Lions, there
are three possibilities for the value of $\gamma:$
\begin{itemize}
\item[$(a)$] Case $1:($\textit{Vanishing}) $\gamma=0.$ Since $P(\omega)$ is non-negative
and nondecreasing, this case is equivalent to saying%
\begin{equation*}
P(\omega)=\lim_{n\to \infty }P_{n}(\omega)=\lim_{n\to \infty
}\sup_{y\in \mathbb{R}}\int_{y-\omega}^{y+\omega}\rho _{n}(x)\ dx=0,
\end{equation*}
for all $\omega<\infty ,\ $or
\item[$(b)$] Case $2:($\textit{Dichotomy}) $\gamma \in (0,\lambda),\ $or
\item[$(c)$] Case $3:($\textit{Compactness}) $\gamma=\lambda,$ that is, there
exists a sequence $\{y_{n}\}$ of real numbers such that $\rho _{n}(.+y_{n})$ is tight,
namely, for all $\varepsilon >0,$ there exists $\omega<\infty $ such that for all $n \in \mathbb{N},$%
\begin{equation*}
\int_{y_{n}-\omega}^{y_{n}+\omega}\rho _{n}(x)dx\geq \lambda-\varepsilon .
\end{equation*}
\end{itemize}
The goal here is to show that the only possibility is the case of compactness, so that the
minimizing sequence $\{(f_{n},g_{n})\}$ for $M_\lambda$ has a subsequence which, up to
translations in the underlying spatial domain, converges strongly in $X$.
We first rule out the vanishing case.

\smallskip

\noindent \textbf{I. The vanishing case does not occur.} We begin with some preliminary lemmas.

\smallskip

\begin{lem}\label{supestimate}
If $g\in C^1(\mathbb{R}), i\in\mathbb{Z}$ and $U_i=[i-1/2,i+1/2]$, then
\[\sup_{r\in U_i}|g(r)|\leq \int_{U_i}[g'_n(s)+g_n(s)]\ ds.\]
\end{lem}
\noindent\textbf{Proof.} Given any $r\in U_i$ and $s\in U_i$, we have
\[g(r)=g(s)+\int_s^r g'(x)\ dx.\]
This implies, for all $i\in\mathbb{Z}$, that
\[|g(r)|\leq |g(s)|+\int_{U_i}|g'(x)|\ dx=|g(s)|+\int_{U_i}|g'(s)|\ ds.\]
Integrating both sides with respect to $s\in U_i$ yields
\[|g(r)||U_i| \leq \int_{U_i} |g(s)|+\int_{U_i}|g'(s)|\ ds |U_i|.\]
Now, using $|U_i|=i+1/2-(i-1/2)=1$ and taking supremum over $r\in U_i$, the
result follows.
\qed

\begin{lem}\label{gamma0}
There exists a $\gamma_0\in(0,\gamma]$ such that
\[\lim_{n\rightarrow \infty} P_n(1/2)=\lim_{n\rightarrow \infty} \sup_{y\in\mathbb{R}}\int_{y-1/2}^{y+1/2}\rho _{n}(x)\ dx\geq \gamma_0.\]
\end{lem}
\noindent\textbf{Proof.} Suppose that we have
\[\lim_{n\rightarrow \infty} P_n(1/2)=\lim_{n\rightarrow \infty} \sup_{y\in\mathbb{R}}\int_{y-1/2}^{y+1/2}\rho _{n}(x)\ dx=0.\]
Letting $U_i=[i-1/2,i+1/2], i\in\mathbb{Z}$ and using Lemma~\ref{supestimate} with $g$ replaced by $g_n^2$, we have
\begin{equation*}
\begin{aligned}
(\sup_{x\in U_i}|g_n(x)|)^2&\leq \int_{U_i}[2 |g_n(s)||g'_n(s)|+(g_n(s))^2]\ ds\\
&\leq\int_{U_i}\left[2\left( \frac{1}{2}\left(|g_n(s)|^2+|g'_n(s)|^2\right)\right)+(g_n(s))^2\right]\ ds\\
&\leq  \int_{i-1/2}^{i+1/2}[(g'_n(s))^2+2(g_n(s))^2]\ ds\\
&\leq C\sup_{y\in\mathbb{R}}\int_{y-1/2}^{y+1/2}\rho _{n}(x)\ dx.
\end{aligned}
\end{equation*}
Hence, we obtain
\begin{equation*}
\begin{aligned}
&0<\lambda =\left|\int_{\mathbb{R}}f_n^2(x)g_n(x)\ dx\right|=\left|\sum_{i=-\infty}^{\infty}\int_{U_i}f_n^2(x)g_n(x)\ dx\right|\\
& \leq \sum_{i=-\infty}^{\infty}\sup_{x\in U_i} |g_n| \int_{U_i}f_n^2(x)\ dx\\
&\leq \left(C\sup_{y\in\mathbb{R}}\int_{y-1/2}^{y+1/2}\rho _{n}(x)\ dx\right)^{1/2}\int_{\mathbb{R}}f_n^2(x)\ dx\\
&\leq \|f_n\|^2 \left(C \sup_{y\in\mathbb{R}}\int_{y-1/2}^{y+1/2}\rho _{n}(x)\ dx\right)^{1/2} \rightarrow 0,
\end{aligned}
\end{equation*}
as $n\rightarrow \infty$, which is a contradiction. Finally,
\[\gamma=\lim_{\omega\rightarrow\infty}P(\omega)\geq P(1/2)=\lim_{n\rightarrow\infty}P_n(1/2)\geq\gamma_0>0.\]
\qed

We now rule out the case of vanishing.

\smallskip

\noindent \textbf{II. The dichotomy case does not occur.} As in the case of vanishing, we first prove some technical lemmas.

\begin{prop}\label{ratio}
For any $\sigma\in\mathbb{R}, 0<M_\sigma<\infty$. Moreover, for any $\sigma_1, \sigma_2>0, M_{\sigma_1}=\left(\frac{\sigma_1}{\sigma_2}\right)^{2/3}M_{\sigma_2}$.
\end{prop}
\noindent\textbf{Proof.}
Clearly, $0\leq I(f,g)<\infty$ for $(f,g)\in X$. Now, $I(f,g)>0$ if either $f\neq 0$ or $g\neq 0$. So, the proposition will hold true after we have shown that given any minimizing sequence
$(f_n,g_n)$ for $M_\sigma$, the sequences $\{f_n\}$ and $\{g_n\}$ cannot vanish. That is,
there exist constants $\gamma_g, \gamma_f> 0$ such that
\[\lim_{n\rightarrow \infty} \|f_n\|\geq \gamma_f \text{ and } \lim_{n\rightarrow \infty} \|g_n\|\geq \gamma_g.\]
In order to prove this statment, we suppose $\displaystyle\lim_{n\rightarrow \infty} \|f_n\|=0$, let $U_i=[i-1/2,i+1/2], i\in\mathbb{Z}$, and use the same argument as in Lemma~\ref{gamma0} to get a contradiction. In the proof of Lemma~\ref{gamma0}, we obtained
\[0<\lambda \leq \|f_n\|^2 \left(C \sup_{y\in\mathbb{R}}\int_{y-1/2}^{y+1/2}\rho _{n}(x)\ dx\right)^{1/2}\rightarrow 0,\]
as $n\rightarrow \infty$, which is a contradiction.
Next, we suppose $\displaystyle\lim_{n\rightarrow \infty} \|g_n\|=0$ and argue similarly. Using the Cauchy-Schwarz inequality, we obtain
\begin{equation*}
\begin{aligned}
&0<\lambda =\left|\int_{\mathbb{R}}f_n^2(x)g_n(x)\ dx\right|=\left|\sum_{i=-\infty}^{\infty}\int_{U_i}f_n^2(x)g_n(x)\ dx\right|\\
& \leq \left|\sum_{i=-\infty}^{\infty}\left(\int_{U_i}(f_n(x))^4\ dx\right)^{1/2}\left(\int_{U_i}(g_n(x))^2\ dx\right)^{1/2}\right|\\
&\leq \sum_{i=-\infty}^{\infty}\left(\sup_{x\in U_i} |f_n(x)|\right)^2 \left(\int_{U_i}(g_n(x))^2\ dx\right)^{1/2}\\
&\leq \left(C\sup_{y\in\mathbb{R}}\int_{y-1/2}^{y+1/2}\rho _{n}(x)\ dx\right)\left(\int_{\mathbb{R}}(g_n(x))^2\ dx\right)^{1/2}\\
&\leq \|g_n\| \left(C \sup_{y\in\mathbb{R}}\int_{y-1/2}^{y+1/2}\rho _{n}(x)\ dx\right) \rightarrow 0,
\end{aligned}
\end{equation*}
as $n\rightarrow \infty$, which is a contradiction. Thus, we have $M_\sigma>0$ for any $\sigma\in\mathbb{R}$. Finally, let $\sigma_1, \sigma_2>0$ and
set $\xi=\left(\frac{\sigma_1}{\sigma_2}\right)^{1/3}$. Clearly, by homegeneity of the functionals, $M_{\sigma_2}=\frac{1}{\xi^2} M_{\xi^3\sigma_2}$, which completes the proof.
\qed

In order to state our next lemma, we will define some new functions. For any $\epsilon>0$, we will first choose a large $\omega\in\mathbb{R}$
and then a large $n\in\mathbb{N}$,
and construct the functions $\eta_{n,\omega}^{(i)}, u_{n,\omega}^{(i)}, i=1,2$ as follows:

Given any $\epsilon>0$, we first find $\omega=\omega(\epsilon)\in\mathbb{R}$ large enough such that $\frac{1}{\omega}\leq \epsilon$ and
\[\gamma-\epsilon/2<P(\omega)\leq P(2\omega)\leq \gamma.\]
Next, we choose $N\in\mathbb{N}$ large enough such that, for all $n\geq N$, we have
\[\gamma-(2\epsilon/3)<P_n(\omega)\leq P_n(2\omega)\leq \gamma +(2\epsilon/3).\]
For each $n\geq N$, we can find $y_n$ such that
\[\int_{y_n-\omega}^{y_n+\omega}\rho _{n}(x)\ dx>\gamma-\epsilon\quad \text{ and } \quad\int_{y_n-\omega}^{y_n+\omega}\rho _{n}(x)\ dx<\gamma+\epsilon.\]
Let $\phi\in C_0^\infty [-2,2]$ and $\psi\in C^\infty (\mathbb{R})$ be such that $\phi\equiv 1$ on $[-1,1]$ and $\phi^2+\psi^2=1$ on $\mathbb{R}$.
Set \[\phi_\omega(x)=\phi\left(\frac{x}{\omega}\right),\quad \psi_\omega(x)=\psi\left(\frac{x}{\omega}\right),\]
and
\[\tilde{\phi}_{n,\omega}(x)=\phi_\omega(x-y_n),\quad
\tilde{\psi}_{n,\omega}(x)=\psi_\omega(x-y_n).\]
Finally, we define
\begin{equation*}
\begin{aligned}
&\eta_{n,\omega}^{(1)}(x)=\tilde{\phi}_{n,\omega}(x)g_n(x), \quad\eta_{n,\omega}^{(2)}(x)=\tilde{\psi}_{n,\omega}(x)g_n(x),\\
& u_{n,\omega}^{(1)}(x)=\tilde{\phi}_{n,\omega}(x)f_n(x), \quad u_{n,\omega}^{(2)}(x)=\tilde{\psi}_{n,\omega}(x)f_n(x),
\end{aligned}
\end{equation*}
and
\[\rho_n^{(i)}(x)=\left(u_{n,\omega}^{(i)'}(x)\right)^2 + \left(\eta_{n,\omega}^{(i)'}(x)\right)^2+\left(u_{n,\omega}^{(i)}(x)\right)^2+\left(\eta_{n,\omega}^{(i)}(x)\right)^2,
i=1,2.\]
\begin{lem}\label{nomegaexist}
For every $\epsilon>0$, there exist a large $\Omega\in\mathbb{R}$ and a large $N\in\mathbb{N}$ such that, for all $n\geq N, \omega\geq\Omega$, we have
\begin{enumerate}[(i)]
\item $I(f_n,g_n)=I\left(u_{n,\omega}^{(1)}, \eta_{n,\omega}^{(1)}\right)+I\left(u_{n,\omega}^{(2)}, \eta_{n,\omega}^{(2)}\right)+O(\epsilon),$
\item $\int_{\mathbb{R}} f_n^2 g_n = \int_{\mathbb{R}} \left(u_{n,\omega}^{(1)}\right)^2 \eta_{n,\omega}^{(1)}
+\int_{\mathbb{R}} \left(u_{n,\omega}^{(2)}\right)^2 \eta_{n,\omega}^{(2)}+O(\epsilon).$
\end{enumerate}
\end{lem}
\noindent\textbf{Proof.} Spelling out and using the Cauchy-Schwarz inequality, we get

\begin{equation*}
\begin{aligned}
&\quad\quad I\left(u_{n,\omega}^{(1)}, \eta_{n,\omega}^{(1)}\right)\\
&=\int_{\mathbb{R}}\left[\left(u_{n,\omega}^{(1)'}\right)^2 + \left(\eta_{n,\omega}^{(1)'}\right)^2
+\alpha_0\left(u_{n,\omega}^{(1)}\right)^2+\beta_0\left(\eta_{n,\omega}^{(1)}\right)^2\right]\\
&=\int_{\mathbb{R}}\bigg[\tilde{\phi}_{n,\omega}^2(g'_n)^2+(\tilde{\phi}'_{n,\omega})^2g_n^2+2\tilde{\phi}_{n,\omega}g'_n\tilde{\phi}'_{n,\omega}g_n
+\tilde{\phi}_{n,\omega}^2(f'_n)^2+(\tilde{\phi}'_{n,\omega})^2f_n^2\\
&\quad\quad\quad +2\tilde{\phi}_{n,\omega}f'_n\tilde{\phi}'_{n,\omega}f_n
+\alpha_0\tilde{\phi}_{n,\omega}^2g_n^2+\beta_0\tilde{\phi}_{n,\omega}^2f_n^2\bigg]\\
&\leq \|\tilde{\phi}_{n,\omega}g'_n\|^2+\|\tilde{\phi}'_{n,\omega}g_n\|^2+2\|\tilde{\phi}_{n,\omega}g'_n\|\|\tilde{\phi}'_{n,\omega}g_n\|
+\|\tilde{\phi}_{n,\omega}f'_n\|^2\\
&\quad +\|\tilde{\phi}'_{n,\omega}f_n\|^2
+2\|\tilde{\phi}_{n,\omega}f'_n\|\|\tilde{\phi}'_{n,\omega}f_n\|
+\alpha_0\|\tilde{\phi}_{n,\omega}g_n\|^2+\beta_0\|\tilde{\phi}_{n,\omega}f_n\|^2.
\end{aligned}
\end{equation*}
Similarly, spelling out $I\left(u_{n,\omega}^{(2)}, \eta_{n,\omega}^{(2)}\right)$ and using $\tilde{\phi}_{n,\omega}^2+\tilde{\psi}_{n,\omega}^2\equiv 1$, we get
\begin{equation*}
\begin{aligned}
&\quad I\left(u_{n,\omega}^{(1)}, \eta_{n,\omega}^{(1)}\right)+I\left(u_{n,\omega}^{(2)}, \eta_{n,\omega}^{(2)}\right)\\
&\leq I(f_n,g_n)+\|\tilde{\phi}'_{n,\omega}g_n\|^2+2\|\tilde{\phi}_{n,\omega}g'_n\|\|\tilde{\phi}'_{n,\omega}g_n\|
+\|\tilde{\phi}'_{n,\omega}f_n\|^2\\
&\quad +2\|\tilde{\phi}_{n,\omega}f'_n\|\|\tilde{\phi}'_{n,\omega}f_n\|+\|\tilde{\psi}'_{n,\omega}g_n\|^2
+2\|\tilde{\psi}_{n,\omega}g'_n\|\|\tilde{\psi}'_{n,\omega}g_n\|\\
&\quad +\|\tilde{\psi}'_{n,\omega}f_n\|^2 +2\|\tilde{\psi}_{n,\omega}f'_n\|\|\tilde{\psi}'_{n,\omega}f_n\|.
\end{aligned}
\end{equation*}
Next, using
\[|\tilde{\phi}'_{n,\omega}|_\infty, |\tilde{\psi}'_{n,\omega}|_\infty\sim O(1/\omega)=O(\epsilon),\]
and, since  $(f_n,g_n)\in X$, that
\[\|f_n\|, \|f'_n\|, \|g_n\|, \|g'_n\|\leq C,\]
for some $C$ independent of $n$, part $(i)$ follows.
For part $(ii)$, we compute
\begin{equation*}
\begin{aligned}
\int_{\mathbb{R}} \left(u_{n,\omega}^{(1)}\right)^2 \eta_{n,\omega}^{(1)}
+\int_{\mathbb{R}} \left(u_{n,\omega}^{(2)}\right)^2 \eta_{n,\omega}^{(2)}=\int_{\mathbb{R}} \left[f_n^2g_n+\left(\tilde{\phi}_{n,\omega}^3+\tilde{\psi}_{n,\omega}^3-1\right)f_n^2g_n\right].
\end{aligned}
\end{equation*}
Upon using $\tilde{\phi}_{n,\omega}\equiv 1, \tilde{\psi}_{n,\omega}\equiv 0$ for $|x-y_n|\leq \omega$ and
$\tilde{\phi}_{n,\omega}\equiv 0, \tilde{\psi}_{n,\omega}\equiv 1$ for $|x-y_n|\geq 2\omega$,
along with the fact that $|\tilde{\phi}_{n,\omega}^3+\tilde{\psi}_{n,\omega}^3-1|\leq 2$, we get
\begin{equation*}
\begin{aligned}
&\left|\int_{\mathbb{R}}\left(\tilde{\phi}_{n,\omega}^3+\tilde{\psi}_{n,\omega}^3-1\right)f_n^2g_n\right|
=\left|\int_{\omega\leq|x-y_n|\leq 2\omega} \left(\tilde{\phi}_{n,\omega}^3+\tilde{\psi}_{n,\omega}^3-1\right)f_n^2g_n\right|\\
&\leq 2|g_n|_\infty\left(\int_{\omega\leq|x-y_n|\leq 2\omega}\rho_n\right)\leq C\epsilon,
\end{aligned}
\end{equation*}
for some $C$ independent of $\omega$ and $n$, and part $(ii)$ is established.\qed

\begin{prop}\label{nondichotomy}
The case $\gamma\in(0,\lambda)$ cannot occur.
\end{prop}
\noindent\textbf{Proof.} Suppose the case $\gamma\in(0,\lambda)$ does occur. Consider a
minimizing sequence $(f_n,g_n)$ for $M_\lambda$. Then, by  Lemma~\ref{nomegaexist}, given $\epsilon>0$,
there exist sequences $\left( u_{n,\omega}^{(i)},\eta_{n,\omega}^{(i)}\right), i=1,2$
such that
\[I(f_n,g_n)=I\left(u_{n,\omega}^{(1)}, \eta_{n,\omega}^{(1)}\right)+I\left(u_{n,\omega}^{(2)}, \eta_{n,\omega}^{(2)}\right)+O(\epsilon).\]
\[\int_{\mathbb{R}} f_n^2 g_n = \int_{\mathbb{R}} \left(u_{n,\omega}^{(1)}\right)^2 \eta_{n,\omega}^{(1)}
+\int_{\mathbb{R}} \left(u_{n,\omega}^{(2)}\right)^2 \eta_{n,\omega}^{(2)}+O(\epsilon).\]
Then, since $\{(f_n,g_n)\}$ is bounded uniformly in $X$, so are  $\left( u_{n,\omega}^{(i)},\eta_{n,\omega}^{(i)}\right), i=1,2$.
Hence, by passing to subsequences but  retaining the same notation, the following limits are well-defined:
\[\sigma_1(\epsilon,\omega)=\lim_{n\to \infty}\int_{\mathbb{R}} \left(u_{n,\omega}^{(1)}\right)^2 \eta_{n,\omega}^{(1)}\text{ and }
\sigma_2(\epsilon,\omega)=\lim_{n\to \infty}\int_{\mathbb{R}} \left(u_{n,\omega}^{(2)}\right)^2 \eta_{n,\omega}^{(2)}.\]
Next, choose a sequence $\epsilon_j\to 0$, and then a sequence $\{\omega_j\}$ depending only on $\{\epsilon_j\}$ so that $\omega_j\to\infty$ , and finally, define
\[\lim_{j\to\infty}\sigma_1(\epsilon_j,\omega_j)=\sigma_1\text{ and } \lim_{j\to\infty}\sigma_2(\epsilon_j,\omega_j)=\sigma_2.\]
Since, $\displaystyle\lim_{n\to\infty}\int_{\mathbb{R}} f_n^2 g_n =\lambda$, we must have $\sigma_1+\sigma_2=\lambda$.
Without loss of generality, suppose $\sigma_2\geq\sigma_1$ as it is only a matter of interchanging $\sigma_1$ and $\sigma_2$.
The following three cases arise:\\

{\textbf Case 1:} $\sigma_1, \sigma_2\in (0,\lambda)$.\\
Taking the limit  $n\to\infty$, and using Proposition~\ref{ratio},
\[I(f_n,g_n)=I\left(u_{n,\omega}^{(1)}, \eta_{n,\omega}^{(1)}\right)+I\left(u_{n,\omega}^{(2)}, \eta_{n,\omega}^{(2)}\right)+O(\epsilon_j)\]
yields
\begin{equation*}
\begin{aligned}
M_\lambda&= M_{\sigma_1(\epsilon_j,\omega_j)}+M_{\sigma_2(\epsilon_j,\omega_j)}+O(\epsilon_j)\\
&=\left[\left(\frac{\sigma_1(\epsilon_j,\omega_j)}{\lambda}\right)^{2/3}+\left(\frac{\sigma_2(\epsilon_j,\omega_j)}{\lambda}\right)^{2/3}\right]M_\lambda+O(\epsilon_j).
\end{aligned}
\end{equation*}
Next, taking the limit $j\to\infty$, we have
\[M_\lambda=\left[\left(\frac{\sigma_1}{\lambda}\right)^{2/3}+\left(\frac{\sigma_2}{\lambda}\right)^{2/3}\right]M_\lambda>M_\lambda,\]
which is a contradiction.

{\textbf Case 2:} $\sigma_1=0$ and $\sigma_2=\lambda$. \\
\begin{equation*}
\begin{aligned}
&\quad\quad I\left(u_{n,\omega}^{(1)}, \eta_{n,\omega}^{(1)}\right)\\
&=\int_{\mathbb{R}}\left[\left(u_{n,\omega}^{(1)'}\right)^2 + \left(\eta_{n,\omega}^{(1)'}\right)^2
+\alpha_0\left(u_{n,\omega}^{(1)}\right)^2+\beta_0\left(\eta_{n,\omega}^{(1)}\right)^2\right]\\
&\geq \min\{1,\alpha_0,\beta_0\}\int_{\mathbb{R}}\left[\left(u_{n,\omega}^{(1)'}\right)^2 + \left(\eta_{n,\omega}^{(1)'}\right)^2
+\left(u_{n,\omega}^{(1)}\right)^2+\left(\eta_{n,\omega}^{(1)}\right)^2\right]\\
&=C\int_{|x-y_n|\leq 2\omega}\left[(f_{n}'(x))^2 + (g_{n}'(x))^2+f_n^2(x)+g_n^2(x)\right]\ dx +O(\epsilon_j)\\
&\geq C\gamma+O(\epsilon_j).
\end{aligned}
\end{equation*}
Taking the limit  $n\to\infty$,
\begin{equation*}
\begin{aligned}
I(f_n,g_n)&=I\left(u_{n,\omega}^{(1)}, \eta_{n,\omega}^{(1)}\right)+I\left(u_{n,\omega}^{(2)}, \eta_{n,\omega}^{(2)}\right)+O(\epsilon_j)\\
&\geq C\gamma +I\left(u_{n,\omega}^{(2)}, \eta_{n,\omega}^{(2)}\right)+O(\epsilon_j)
\end{aligned}
\end{equation*}
yields
\begin{equation*}
\begin{aligned}
M_\lambda&\geq C\gamma +M_{\sigma_2(\epsilon_j,\omega_j)}+O(\epsilon_j).
\end{aligned}
\end{equation*}
Upon taking the limit  $j\to\infty$,
\[M_\lambda\geq  C\gamma+M_\lambda>M_\lambda,\]
which is a contradiction.

\textbf{ Case 3:} $\sigma_1 <0$ and $\sigma_2=\lambda-\sigma_1>\lambda$. \\
Taking the limit  $n\to\infty$, and using Proposition~\ref{ratio},
\begin{equation*}
\begin{aligned}
I(f_n,g_n)&=I\left(u_{n,\omega}^{(1)}, \eta_{n,\omega}^{(1)}\right)+I\left(u_{n,\omega}^{(2)}, \eta_{n,\omega}^{(2)}\right)+O(\epsilon_j)\\
&\geq I\left(u_{n,\omega}^{(2)}, \eta_{n,\omega}^{(2)}\right)+O(\epsilon_j)
\end{aligned}
\end{equation*}yields
\begin{equation*}
\begin{aligned}
M_\lambda&\geq M_{\sigma_2(\epsilon_j,\omega_j)}+O(\epsilon_j)\\
&=\left(\frac{\sigma_2(\epsilon_j,\omega_j)}{\lambda}\right)^{2/3}M_\lambda+O(\epsilon_j).
\end{aligned}
\end{equation*}
Finally, taking the limit  $j\to\infty$,
\[M_\lambda\geq\left(\frac{\sigma_2}{\lambda}\right)^{2/3}M_\lambda>M_\lambda,\]
which is a contradiction.\qed

\medskip

\noindent \textbf{III. Proof of Theorem~\ref{existence}.}
As we ruled out both vanishing and dichotomy, Lions'
concentration compactness lemma guarantees that the sequence $
\{\rho _{n}\}$ is tight up to translation. Thus, there exists a sequence of real numbers $
\{y_{n}\}$ such that given any $\varepsilon >0,$ we can find $\omega=\omega(\varepsilon
)\geq \frac{1}{\varepsilon}$ such that for all $n\in \mathbb{N},$
\begin{equation*}
\begin{aligned}
\int_{|x-y_n|\leq\omega}\rho _{n}(x)dx
\geq\gamma-\varepsilon ,\; \int_{|x-y_n|\geq\omega}\rho _{n}(x)\ dx\leq\varepsilon,
\end{aligned}
\end{equation*}

and
\begin{equation}\label{outsideball}
\begin{aligned}
\left|\int_{|x-y_n|\geq\omega}f_n^2(x)g_n(x)\ dx\right|
&\leq C \|(f_n,g_n)\|_X \int_{|x-y_n|\geq\omega}\rho _{n}(x)dx\\
&=O(\varepsilon).
\end{aligned}
\end{equation}
Denote $\tilde{f_n}=f_n(\cdot+y_n)$ and $\tilde{g_n}=g_n(\cdot+y_n)$. Then, by \eqref{outsideball},
\begin{equation*}
\begin{aligned}
\left|\int_{-\omega}^{\omega}\tilde{f_n}^2\tilde{g_n}-\lambda\right|
&=\left|\int_{|x-y_n|\geq\omega}f_n^2(x)g_n(x)\ dx-\int_{\mathbb{R}}f_n^2(x)g_n(x)\ dx\right|\\
&=\left|\int_{|x-y_n|\leq\omega}f_n^2(x)g_n(x)\ dx\right|=O(\varepsilon).
\end{aligned}
\end{equation*}
Thus,
\begin{equation}\label{squeeze}
\lambda-O(\varepsilon)\leq\int_{-\omega}^{\omega}\tilde{f_n}^2\tilde{g_n}\leq\lambda+O(\varepsilon).
\end{equation}
Since $\{(\tilde{f_n},\tilde{g_n})\}$
is uniformly bounded in $X$, by Banach-Alaoglu's theorem, there exists a subsequence, again labeled
$\{(\tilde{f_n},\tilde{g_n})\}$ for notational convenience, that converges weakly in $X$, say, to $(\phi, \psi)$.
Using the Cauchy-Schwarz inequality and the compact embedding of $H^1([-\omega,\omega])$ in $L^2([-\omega,\omega])$, we have
\begin{equation*}
\begin{aligned}
&\quad\int_{-\omega}^{\omega}|\tilde{f_n}^2\tilde{g_n}-\phi^2\psi|
=\int_{-\omega}^{\omega}|(\tilde{f_n}+\phi)\tilde{g_n}(\tilde{f_n}-\phi)+\phi^2(\tilde{g_n}-\psi)|\\
&\leq |\tilde{f_n}+\phi|_\infty \cdot \|\tilde{g_n}\|\cdot |\tilde{f_n}-\phi|_{L^2(-\omega,\omega)}
+\|\phi\|^2\cdot |\tilde{g_n}-\psi|_{L^2(-\omega,\omega)}\\
&\leq C \left(|\tilde{f_n}-\phi|_{L^2(-\omega,\omega)}+|\tilde{g_n}-\psi|_{L^2(-\omega,\omega)}\right)\to 0,
\end{aligned}
\end{equation*}
as $n\to\infty$. Hence, there exists a large enough $N$ such that for all $n\geq N$, \eqref{squeeze} yields
\begin{equation*}
\lambda-O(\varepsilon)\leq\int_{-\omega}^{\omega}\phi^2\psi\leq\lambda+O(\varepsilon).
\end{equation*}
Now, choosing $\varepsilon_j=\frac{1}{j}, j\in\mathbb{N}$ and determining $\omega_j=\omega_j(\varepsilon_j)\geq \frac{1}{\varepsilon_j}=j$ and letting $j\to\infty$, we obtain
\[\int_{\mathbb{R}}\phi^2\psi=\lambda.\]
Moreover, since $I$ is weak lower semi-continuous and invariant under translation,
\[M_\lambda=\lim_{n\to\infty} I(\tilde{f_n},\tilde{g_n})\geq I(\phi,\psi)\geq M_\lambda.\]
This means $(\phi,\psi)\in \mathcal{S}_\lambda$. Since $(\phi,\psi)$ is a solution of the variational
problem of minimizing $I(f,g)$ subject to the constraint $\int_{\mathbb{R}} f^2g=\lambda$, it must satisfy the corresponding Euler-Lagrange equation. That is, there exists some Lagrange multiplier $\kappa\in\mathbb{R}$ such that
\begin{equation}\label{eulerlag}
\left\{
\begin{aligned}
& -\phi^{\prime \prime}+\alpha_0 \phi =\kappa\phi \psi,\\
& -\psi^{\prime \prime}+\beta_0 \psi =\frac{\kappa}{2}\phi^2.
\end{aligned}
\right.
\end{equation}
Multiplying the first and second equations by $\phi$ and $\psi$ respectively and integrating over $\mathbb{R}$ yields
\begin{equation*}
\left\{
\begin{aligned}
& \int_{\mathbb{R}}(\phi^{\prime})^2+\alpha_0 \phi^2 =\kappa\int_{\mathbb{R}}\phi^2 \psi,\\
& \int_{\mathbb{R}}(\psi^{\prime})^2+\beta_0 \psi^2 =\frac{\kappa}{2}\int_{\mathbb{R}}\phi^2\psi.
\end{aligned}
\right.
\end{equation*}
Adding the two equations together and using the fact that $(\phi,\psi)\in \mathcal{S}_\lambda$, we get
\[\kappa=\frac{2M_\lambda}{3\lambda}>0.\]
Clearly, $(\kappa\phi,\kappa\psi)$ solves the system \eqref{auxsys}, or, equivalently,
\[W(x,t)=\kappa e^{i\sigma t}\phi(x), \:  V(x,t)=\kappa e^{2i\sigma t}\psi(x)\]
are solitary wave solutions to the system \eqref{nlseq}.
Finally, we establish the smoothness of $\phi$ and $\psi$ using a standard bootstrap argument.
The system \eqref{auxsys} can be rewritten as
\begin{equation}\label{convol}
\left\{
\begin{aligned}
&\phi= K_{\alpha_0}\star\phi\psi\,\\
& \psi=\frac{1}{2} \;K_{\beta_0}\star\phi^2,
\end{aligned}
\right.
\end{equation}
where given any $s>0$, the kernel $K_s$ is defined via its Fourier transform as
\[q(\zeta)=\hat{K_s}(\zeta)=\frac{1}{s+\zeta^2}.\]
Since $H^s(\mathbb{R}^n)$ is an algebra if $s>n/2$, we have $\phi\psi,\psi^2\in H^1(\mathbb{R})$. Since the convolution operation takes $H^s$ to $H^{s+2}$ for
any $s\geq 0$, \eqref{convol} implies that $\phi,\psi\in H^3$. Again, using the fact that $\phi\psi,\psi^2\in H^3$, we obtain $\phi,\psi\in H^5$. Continuation of this induction argument leads to $\phi,\psi\in H^\infty$. Exponetial decay of $\phi$ and $\psi$ at infinity follows from another standard technique. Using Theorem 8.1.1 in
\cite{[Caz]}, there exits a $\delta>0$ such that
$e^{\delta|\cdot|}\phi, e^{\delta|\cdot|}\psi\in L^{\infty}$. Thus,
\[|\phi(x)|\leq \left|e^{\delta|\cdot|}\phi\right|_\infty e^{-\delta|x|}\to 0,\; |\psi(x)|\leq \left|e^{\delta|\cdot|}\psi\right|_\infty e^{-\delta|x|}\to 0,\]
as $|x|\to\infty$.
\qed

In the course of the proof of Theorem~\ref{existence} above, we have also managed to establish the following result:

\begin{prop}
If $\lambda>0$, the Lagrange multiplier $\kappa\in\mathbb{R}$ associated to the variational
problem of minimizing $I(f,g)$ subject to the constraint $\int_{\mathbb{R}} f^2g=\lambda$ is positive.
\end{prop}

\medskip


\section*{Acknowledgments}
The author is grateful to JCHS Faculty
Development Grant Program for a grant that partially supported this project.

\end{document}